\documentclass[amsthm]{elsart}
\usepackage{yjsco}
\usepackage{natbib}
\usepackage{amsmath,amssymb}
\usepackage{graphicx}
\usepackage{multirow}
\usepackage{enumerate}
\usepackage{array}
\usepackage{algorithm}
\usepackage{algorithmic}
\usepackage{url}
\usepackage{hyperref}
\def\bibsection{\section*{References}}

\begin{document}

\begin{frontmatter}

\title{Gr\"obner bases and cocyclic Hadamard matrices}

\author{V. \'Alvarez}\ead{valvarez@us.es}\address{}
\vspace{-1.25cm}

\author{J.A. Armario}\ead{armario@us.es}\address{}

\vspace{-1.25cm}
\author{R.M. Falc\'{o}n}\ead{rafalgan@us.es}\address{}

\vspace{-1.25cm}
\author{M.D. Frau}\ead{mdfrau@us.es}\address{}

\vspace{-1.25cm}
\author{F. Gudiel}\ead{gudiel@us.es}
\address{Dpto. Matem\'{a}tica Aplicada I. Univ. Sevilla.\\Avda. Reina Mercedes s/n 41012 Sevilla, Spain}

\begin{abstract}
Hadamard ideals were introduced in 2006 as a set of nonlinear polynomial equations whose zeros are uniquely related to Hadamard matrices with one or two circulant cores of a given order. Based on this idea, the cocyclic Hadamard test enable us to describe a polynomial ideal that characterizes the set of cocyclic Hadamard matrices over a fixed finite group $G$ of order $4t$. Nevertheless, the complexity of the computation of the reduced Gr\"obner basis of this ideal is $2^{O(t^2)}$, which is excessive even for very small orders. In order to improve the efficiency of this polynomial method, we take advantage of some recent results on the inner structure of a cocyclic matrix to describe an alternative polynomial ideal that also characterizes the mentioned set of cocyclic Hadamard matrices over $G$. The complexity of the computation decreases in this way to $2^{O(n)}$, where $n$ is the number of $G$-coboundaries. Particularly, we design two specific procedures for looking for $\mathbb{Z}_t \times \mathbb{Z}_2^2$-cocyclic Hadamard matrices and $D_{4t}$-cocyclic Hadamard matrices, so that larger cocyclic Hadamard matrices (up to $t \leq 31$) are explicitly obtained.
\end{abstract}

\end{frontmatter}

\newtheorem{lemma}[thm]{Lemma}

\newcommand{\D}{\displaystyle}

\section{Introduction.}

A {\it Hadamard matrix} $H$ of order $n$ is an $n\times n$ matrix with every entry either $1$ or $-1$, which satisfies $HH^T=nI$, where $I$ is the identity matrix of order $n$. Although it is well-known that $n$ has to be necessarily $1,\,2$ or a multiple of $4$ (as soon as three or more rows have to be simultaneously orthogonal), there is no certainty whether such a Hadamard matrix exists at every possible order. Currently, the smallest order for which no Hadamard matrix is known is 668. The {\it Hadamard conjecture} asserts that there exists a Hadamard matrix of order $4t$ for every natural number $t$.

\vspace{0.2cm}

There exist so many different constructions for Hadamard matrices: Sylvester, Paley, Williamson, Ito, Goethals-Seidel, one and two circulant cores or cocyclic matrices, amongst others (see \cite{Hor07}). Nevertheless, most of them fail to yield Hadamard matrices for every order which is a multiple of 4 and therefore are not suitable candidates for a proof of the Hadamard conjecture. Among all these constructions, it seems that the two most promising are the Goethals-Seidel arrays and the cocyclic constructions. Actually, the one and two circulant cores constructions have recently been described to be somehow cocyclic-based (the cores themselves are cocyclic over $\mathbb{Z}_{4t-1}$ and $D_{4t-2}$, respectively). A stronger version of the Hadamard conjecture, posed by \cite{HdL95} is the {\em cocyclic Hadamard conjecture:} this states that there exists a cocyclic Hadamard matrix at every possible order. Currently the smallest order for which no cocyclic Hadamard matrix is known is 188.

\vspace{0.2cm}

\cite{Kotsireas2006} introduced the concept of {\em Hadamard ideal} as a set of nonlinear polynomial equations whose zeros determine the set of Hadamard matrices with one  circulant core. Shortly after, \cite{Kotsireas2006a} used the same ideal together with a series of new polynomials in order to determine the set of Hadamard matrices with two circulant cores, by means of which they computed the Hadamard matrices with two circulant cores up to order 52.

\vspace{0.2cm}

In this paper, we define several {\em cocyclic Hadamard ideals}, whose zeros determine the set of cocyclic Hadamard matrices over a finite group $G$ of order $4t$. Based on the cocyclic test of \cite{HdL95}, our first approach (Theorem \ref{thm1}) gives rise to a procedure {\tt CocGM($t,G,opt$)} which  works just for very small $t$, actually $t \leq 3$. In order to improve the efficiency of this polynomial method and provided a basis of $G$-cocycles is known, we define in Theorem \ref{thmbasesgeneral} an alternative ideal based upon the system of equations described by \cite{AAFR08}, which also characterize the set of $G$-cocyclic Hadamard matrices. This gives a procedure {\tt CocCB($t,G,opt$)} suitable for larger values of $t$. Furthermore, from the knowledge of the properties of cocyclic matrices over $\mathbb{Z}_t\times \mathbb{Z}_2^2$ and $D_{4t}$ described by \cite{AGG15, AAFGGO16}, improved versions of this procedure ({\tt CocAH($t,col,dist,H$)} and {\tt CocDH($t,opt,H$)}, based on Theorems \ref{thm2} and \ref{thm3}, respectively) are used to perform local searches for $\mathbb{Z}_t\times \mathbb{Z}_2^2$-cocyclic Hadamard matrices and $D_{4t}$-cocyclic Hadamard matrices, so that matrices of order up to $4t \leq 124$ are found.

\vspace{0.2cm}

All the procedures have been implemented as a library {\em hadamard.lib} in the open computer algebra system for polynomial computations {\sc Singular}, developed by \cite{Decker2016}. Examples illustrating the use of this library and the library itself are available online at \texttt{http://personales.us.es/raufalgan/LS/hadamard.lib}.

\vspace{0.2cm}

The remainder of the paper is organized as follows. The first part of Section \ref{sec:prelim} is devoted to describe some preliminary concepts and results on Hadamard matrices and algebraic geometry, that are used in the rest of the paper. Later, we define a zero-dimensional ideal that determines the set of cocyclic Hadamard matrices over a given group $G$ of order $4t$, which comes from a straightforward translation of the cocyclic Hadamard test of \cite{HdL95}. In Section \ref{sec:cocyclic}, we propose an alternative to the previous construction by defining a new zero-dimensional ideal, based on the results of \cite{AAFR08}. Actually, we particularize this procedure for  $\mathbb{Z}_t\times \mathbb{Z}_2^2$-cocyclic Hadamard matrices and $D_{4t}$-cocyclic Hadamard matrices, attending to the properties described by \cite{AGG15, AAFGGO16}. The last section is devoted to conclusions and outlines for further work.

\section{Preliminaries.}\label{sec:prelim}

We expose in this section some basic concepts and results on Hadamard matrices and algebraic geometry that are used throughout the paper. We refer to the monographs of \cite{McL95}, \cite{Hor07}, \cite{dLF11} and \cite{Cox1998, Cox2007} for more details about these topics.

\vspace{0.2cm}

Assume throughout that $G=\{ g_1=1, \ldots, g_{4t}\}$ is a multiplicative finite group of $4t$ elements, not necessarily abelian. A function $\psi\colon G\times G\rightarrow \langle -1\rangle\cong \mathbb{Z}_2$ is said to be a {\em (binary) cocycle} over $G$, or simply $G$-cocycle for short,  if it satisfies that
\begin{equation}\label{eqcocyclic}
\psi(g_i,g_j)\psi(g_ig_j,g_k)=\psi(g_j,g_k)\psi(g_i,g_jg_k), \text{ for all } g_i,g_j,g_k\in G.
\end{equation}

\vspace{0.2cm}

The cocycle $\psi$ is naturally displayed as a {\it cocyclic matrix} $M_\psi$ whose $(i,j)^{\text{th}}$ entry is $\psi(g_i,g_j)$ for all $g_i,g_j\in G$. Since it must be $\psi(1,g_j)=\psi(g_i,1)$ for all $g_i,g_j\in G$, the first row and column of $M_\psi$ are all either 1 or $-1$. In the first case, the cocycle $\psi$ and its cocyclic matrix $M_\psi$ are said to be {\it normalized}. There is a one to one correspondence between normalized and non normalized cocycles. In what follows, we reduce ourselves to normalized cocycles, for commodity.

\vspace{0.2cm}

Given $g_d\in G$, the {\em elementary coboundary} $\partial_d$ is the cocycle over $G$ defined as
\[\partial_d(i,j):=\delta_{g_d}(g_i)\delta_{g_d}(g_j)\delta_{g_d}(g_ig_j),\]
where $\delta_{g_d}\colon G\rightarrow  \{-1,1\}$ is the characteristic set map such that $\delta_{g_d}(g_i)=-1$ if $g_i=g_d$ and $1$, otherwise. The {\it generalized coboundary matrix} $\overline{M}_{\partial_d}$ consists of negating the $d^{\text{th}}$-row of the matrix $M_{\partial_d}$. Note that negating a row or a column of a matrix does not change its Hadamard character. This is just a particular case of a more general set: there is an equivalence relation (termed {\em Hadamard equivalence}) on Hadamard matrices, so that two matrices are Hadamard equivalent whenever they differ in a series of row and/or column negations and/or permutations.

\vspace{0.2cm}

\cite{AAFR08} proved that every generalized coboundary matrix $\overline{M}_{\partial_d}$ has the following properties

\vspace{0.2cm}

\begin{enumerate}[a)]
\item $\overline{M}_{\partial_d}$ contains exactly two negative entries in each row $s\neq 1$, which are located at positions $(s,d)$ and $(s,e)$, for $g_e=g_s^{-1}g_d$.
\item Given $g_s\neq 1$ and $g_c$ in $G$, there are exactly two generalized coboundary matrices ($\overline{M}_{\partial_c}$ and $\overline{M}_{\partial_d}$), with a negative entry in the position $(s,c)$, where $g_d=g_sg_c$.
\item Two generalized coboundary matrices share their two negative entries at the $s^{\text{th}}$ row if and only if $g_s^2=1$.
\end{enumerate}

\vspace{0.2cm}

A {\em basis} ${\bf B}=\{ \psi_1 , \ldots, \psi_k\}$ of cocycles over $G$ consists of some elementary coboundaries $\partial_i$ and some representative cocycles. Every cocycle over $G$ admits a unique representation as a product of the generators in ${\bf B}$, $\psi=\psi_1^{x_1}\cdots \psi_k^{x_k}$, $x_i \in \mathbb{Z}_2$. The tuple $(x_1,\ldots,x_k)_{\bf B}$ defines the {\em coordinates} of $\psi$ with regards to ${\bf B}$. Accordingly, every cocyclic matrix $M_\psi=(\psi(i,j))$, for $\psi=(x_1,\ldots,x_k)_{\bf B}$, admits a unique decomposition $M_\psi=M_{\psi_1}^{x_1} \cdots  M_{\psi_k}^{x_k}$ as the Hadamard pointwise product of those matrices $M_{\psi_i}$ corresponding to entries $x_i=1$. In general, the number of elements of ${\bf B}$ is an open question. Furthermore, since the elementary coboundary $\partial_1$ related to the identity element $1\in G$ is not normalized, we may assume that $\partial_1 \notin {\bf B}$. In what follows, we use generalized coboundary matrices instead of classical coboundary matrices. Let us point out that any matrix obtained as the Hadamard product of generalized coboundary matrices and representative cocycles is Hadamard equivalent to a cocyclic matrix by means of negations of certain rows.

\vspace{0.2cm}

A cocycle $\psi$ (over $G$) is said to be {\em orthogonal} if its cocyclic matrix $M_\psi$ is Hadamard. In such a case, $M_\psi$ is said to be a {\it a cocyclic Hadamard matrix over $G$}. The set of cocyclic Hadamard matrices over $G$ is denoted by $\mathcal{H}_G$. The {\em cocyclic Hadamard test} of \cite{HdL95} asserts that a normalized cocyclic matrix $M_\psi$ is Hadamard if and only if
\begin{equation}\label{cTestHad}
\sum_{j\in G} \psi(i,j)=0, \text{ for all } i\in G\setminus\{1\}.
\end{equation}
A row of $M_\psi$ is termed {\em Hadamard row} precisely when its summation is zero. This way, $M_\psi$ is Hadamard if and only if every row (but the first) is a Hadamard row.

\vspace{0.2cm}

We expose now some basic concepts of algebraic geometry. Let $X$ and $\mathbb{K}[X]$ be, respectively, the set of $m$ variables $\{x_1,\ldots,x_m\}$ and the related multivariate polynomial ring over a field $\mathbb{K}$. The {\em affine variety} $V(I)$ of an ideal $I\subseteq \mathbb{K}[X]$ is the set of points in $\mathbb{K}^m$ that are zeros of all the polynomials of $I$. The ideal $I$ is said to be {\em zero-dimensional} if $V(I)$ is finite. It is said to be {\em radical} if every polynomial $p\in \mathbb{K}[X]$ belongs to $I$ whenever there exists a natural number $n$ such that $p^n\in I$. A {\em term order} $<$ on the set of monomials of $\mathbb{K}[X]$ is a multiplicative well-ordering that has the constant monomial $1$ as its smallest element. The largest monomial of a polynomial $p$ of $I$ with respect to the term order $<$ is its {\em leading monomial}. The ideal generated by the leading monomials of all the non-zero elements of $I$ is its {\em initial ideal} $I_<$. Those monomials of polynomials of $I$ that are not leading monomials of any polynomial of $I$ are called {\em standard monomials}. If the ideal $I$ is zero-dimensional, then the number of standard monomials of $I$ coincides with the dimension of $\mathbb{K}[X]/I$ over $\mathbb{K}$, which is greater than or equal to the number of points of $V(I)$. The equality holds when $I$ is radical. This dimension can be obtained by computing the {\em Hilbert function} $\mathrm{HF}_{\mathbb{K}[X]/I}$, which maps each non-negative integer $d$ onto $\mathrm{dim}_k(\mathbb{K}[X]_d/I_d)$, where $\mathbb{K}[X]_d$ denotes the set of homogeneous polynomials in $\mathbb{K}[X]$ of degree $d$ and $I_d=\mathbb{K}[X]_d\cap I$. In particular, $\mathrm{dim}_k(\mathbb{K}[ X]/I)=\sum_{0\leq d}\mathrm{HF}_{\mathbb{K}[X]/I}(d)$. If the ideal $I$ is zero-dimensional, then the number $\mathrm{HF}_{\mathbb{K}[X]/I}(d)$ coincides with the set of standard monomials of degree $d$, regardless of the term order. As a consequence, the Hilbert function of $\mathbb{K}[X]/I$ coincides with that of $\mathbb{K}[X]/I_<$, for any term order $<$, which can be obtained by using for instance the algorithm of \cite{Mora1983}. Previously, it is required to determine the initial ideal $I_<$. In any case, \cite{Bayer1992} already proved that the problem of computing Hilbert functions is NP-complete.

\vspace{0.2cm}

A {\em Gr\"obner basis} of the ideal $I$ is any subset $GB$ of polynomials of $I$ whose leading monomials with respect to a given term order generate the initial ideal $I_<$. It is {\em reduced} if all its polynomials are monic and no monomial of a polynomial in $GB$ is generated by the leading monomials of the rest of polynomials in the basis. There exists only one reduced Gr\"obner basis of the ideal $I$. This basis generates the initial ideal $I_<$ and can be used, therefore, to determine the cardinality of its affine variety $V(I)$. Further, the points of this variety can been enumerated once the reduced Gr\"obner basis is decomposed into finitely many disjoint subsets, each of them being formed by the polynomials of a triangular system of polynomial equations, whose factorization and subsequent resolution are easier than the system related to the generators of the original ideal $I$. See in this regard the articles of \cite{Hillebrand1999}, \cite{Lazard1992} and \cite{Moller1993}.

\vspace{0.2cm}

Gr\"obner bases can be used, therefore, to determine both the cardinality and the elements of the set $\mathcal{H}_G$ of cocyclic Hadamard matrices over a multiplicative finite group $G$ of $4t$ elements. Let $\mathbb{Q}[X_G]$ be the polynomial ring over the set of variables $\{X_G\}=\{x_{i,j}\colon\, g_i,g_j\in G\}$ and let us define the polynomial
\[p_{i,j,k}(\text X):=x_{i,j}-x_{ij,k}x_{j,k}x_{i,jk}, \text{ for all } g_i,g_j,g_k\in G,\]
where the products $ij$ and $jk$ are induced by the group law in $G$. The next result shows how the set $\mathcal{H}_G$ of cocyclic Hadamard matrices over $G$ can be identified with the affine variety defined by a zero-dimensional radical ideal of nonlinear polynomials in $\mathbb{Q}[X_G]$.

\begin{thm} \label{thm1} The set $\mathcal{H}_G$ can be identified with the set of zeros of the zero-dimensional ideal
\begin{align*}I_G\, := \, & \langle\,x_{i,j}^2-1, x_{1,j}-1, x_{i,1}-1\colon\, i,j\in G\setminus\{1\}\,\rangle\, + \langle\,p_{i,j,k}(\text X)\colon\, i,j,k\in G\,\rangle\, + \\ \ & \langle\,\sum_{j\in G}x_{i,j}\colon\, i\in G\setminus\{1\}\,\rangle\, \subset \mathbb{Q}[X_G].
\end{align*}

Besides, $|\mathcal{H}_G|= \mathrm{dim}_{\mathbb{Q}}(\mathbb{Q}[X_G]/I_G)$.
\end{thm}

\begin{pf}
Let $P=(p_{i,j})$ be a point of the affine variety $V(I_G)$. From the first subideal of $I_G$, every component $p_{i,j}$ of $P$ is either $1$ or $-1$, for all $i,j\in G$. Let $\psi:G\times G\rightarrow \{\pm 1\}$ be defined such that $\psi(i,j)=p_{i,j}$, for all $g_i,g_j\in G$. Since the second subideal of $I_G$ implies that $\psi$ satisfies identity (\ref{eqcocyclic}) for all $g_i,g_j,g_k\in G$, the point $P$ can be identified with the cocyclic matrix $M_\psi$ related to $\psi$. Finally, the third subideal of $I_G$ implies that $M_\psi$ satisfies identity (\ref{cTestHad}) and hence, $M_\psi$ is Hadamard. The affine variety $V(I_G)$ coincides, therefore, with the set $\mathcal{H}_G$, whose finiteness involves
the ideal $I_G$ to be zero-dimensional. Besides, since $I_G\cap \mathbb{Q}[x_{i,j}]=\langle\,x_{ij}^2-1\,\rangle\subseteq I_G$ for all $i,j\in G$, Proposition 2.7 of \cite{Cox1998} involves $I_G$ to be radical and hence, $|\mathcal{H}_G|=|V(I_G)|= \mathrm{dim}_{\mathbb{Q}}(\mathbb{Q}[X_G]/I_G)$.
\end{pf}

\vspace{0.5cm}

Notice that the computation of the reduced Gr\"obner basis of a zero-dimensional ideal is extremely sensitive to the number of variables. See in this regard the articles of \cite{Hashemi2009}, \cite{Hashemi2011}, \cite{Lakshman1991} and \cite{Lakshman1991a}. In the last reference, the authors proved that the complexity of our computation is $d^{O(n)}$, where $d$ is the maximal degree of the polynomials of the ideal and $n$ is the number of variables. In the case of Theorem \ref{thm1}, this complexity is $2^{O(t^2)}$, which renders the computation only possible for very low values of $t$ that are not useful to analyze the Hadamard conjecture. Thus, for instance, using our procedure {\tt CocGM($t,G,opt$)} in an Intel Core i7-2600, with a 3.4 GHz processor and 16 GB of RAM, the computation of the reduced Gr\"obner bases of the ideals related to the group $\mathbb{Z}_t\times \mathbb{Z}^2_2$ and the dihedral group $D_{4t}$ are only feasible for $t\leq 3$. For higher orders, the system runs out of memory. Notice that, depending on whether the parameter {\tt opt} is equal to 1 or 2, the procedure calculates either just the number of cocyclic Hadamard matrices over $G$ or the explicit full set of these matrices.

\vspace{0.2cm}

In Section \ref{sec:cocyclic} we define another ideal $J_G$ for computing $\mathcal{H}_{G}$ in a more subtle way, progressing on the previous work of \cite{AAFR08} and provided a basis for cocycles over $G$ is known. Unfortunately, it will be still extremely hard to compute $\mathcal{H}_{G}$ for large $|G|$. Nevertheless, taking advantage of the properties of cocyclic matrices over $D_{4t}$ and $\mathbb{Z}_t\times \mathbb{Z}^2_2$ described by \cite{AGG15, AAFGGO16}, this ideal $J_G$ may be specifically simplified for computing $\mathcal{H}_{D_{4t}}$ and $\mathcal{H}_{\mathbb{Z}_t\times \mathbb{Z}^2_2}$ in a better way.

\section{Ideals built from a basis for $G$-cocycles}\label{sec:cocyclic}

In order to reduce the complexity of the computation of the reduced Gr\"obner basis that has been exposed in the previous section, we consider a new zero-dimensional radical ideal $J_G$ related to the set $\mathcal{H}_G$, where we diminish the number of variables and the maximal degree of the polynomials progressing on the knowledge of an explicit basis for cocycles over $G$.

\vspace{0.2cm}

Let $G$ be a multiplicative finite group of order $4t$, ${\bf B}=\{\psi_1,\ldots,\psi_k\}$ be a basis for normalized cocycles over $G$ and $\psi$ be a normalized cocycle over $G$ of coordinates $(x_1,\ldots,x_k)_{\bf B}$ with regards to ${\bf B}$. Let $m^d_{i,j}$ denote the $(i,j)^{\text{th}}$ entry of $M_{\psi_d}$, so that the $(i,j)^{\text{th}}$ entry of $M_\psi$ is $(m_{i,j}^1)^{x_1}\cdots (m_{i,j}^k)^{x_k}$. Recall that normalized cocyclic Hadamard matrices are precisely those matrices that are built up from Hadamard rows (excepting the first row, consisting  all of 1s). In these circumstances, the $i^{th}$-row of the previous matrix $M_\psi$ is Hadamard if and only if
\[\D\sum _{j=1}^{4t}(m_{i,j}^1)^{x_1}\cdots (m_{i,j}^k)^{x_k}=0.\]
The next result holds.

\begin{thm}[\cite{AAFR08}]\label{sistema}
The matrix $M_\psi$ is Hadamard if and only if the vector of coordinates $(x _1, \ldots , x _k)_{\bf B}$ of $\psi$ with regards to ${\bf B}$
satisfies the following system of $4t-1$ equations and $k$ unknowns

\begin{equation} \label{system}\left\{ \begin{array}{lcc} (m_{2,1}^1)^{x_1}\ldots (m_{2,1}^k)^{x_k}+\ldots+
(m_{2,4t}^1)^{x_1} \ldots (m_{2,4t}^k)^{x_k}&=&0 \\ &\vdots& \\
(m_{4t,1}^1)^{x_1} \ldots (m_{4t,1}^k)^{x_k}+\ldots+ (m_{4t,4t}^1)^{x_1} \cdots (m_{4t,4t}^k)^{x_k}&=&0
\end{array}\right.\end{equation} \hfill $\Box$
\end{thm}

\vspace{0.5cm}

The solutions of the system (\ref{system}) constitute precisely the whole set of normalized cocyclic Hadamard matrices over $G$. Trying to solve this system may be as complicated as performing an exhaustive search for cocyclic Hadamard matrices over $G$. Instead, we intend to translate the system (\ref{system}) in terms of a set of nonlinear $\mathbb{Q}[X]$-polynomial equations over the set of variables $\{ X\}=\{x_1,\ldots,x_k\}$ (whose $0,1$ values are related to the coordinates of  $G$-cocycles with regards to ${\bf B}$), and to study the structure of the associated ideal.

\vspace{0.2cm}

A succinct algebraic description of the quadratic constrains  $\{ X\}\subset \{0,1\}^k$ is provided by the following set of $k$ algebraic equations:
\begin{equation}\label{0dime}
x_i(x_i-1)=0, \,\, \text{ for all } i\in\{1,\ldots,k\}.
\end{equation}
In order to define the rest of polynomial equations that arise from the system (\ref{system}), we use the next two main ideas or simplifications:

\vspace{0.2cm}

\begin{itemize}
\item From a practical point of view, we may assume we work with a fixed representative cocycle $\rho$ among all of the possible choices of representative cocycles. In fact, empirically, in the groups most intensively studied, there always exists a better combination of representative cocycles for producing Hadamard matrices. See in this regard the works of \cite{AAFR08, AGG15, AAFGGO16}, \cite{BH95}, \cite{Fla97} and \cite{Hor07}. We will denote by $M_\rho=(r_{i,j})$ the matrix related to the Hadamard product $\rho$ of these representative cocycles. Obviously, this pruning in the searching space might eliminate some cocyclic Hadamard matrices. If we want to find the whole set of cocyclic Hadamard matrices, we have to perform an analogous search for the other possible choices of $M_\rho$. In what follows we assume that $\psi _1, \ldots , \psi _{k-m} \in {\bf B}$ are $G$-coboundaries, $\psi _{k-m+1}, \ldots, \psi_k\in {\bf B}$ are representative $G$-cocycles and $\rho= \prod _{i=k-m+1}^k \psi _i ^{x_i}$ is a fixed linear combination of these representative cocycles.
\item  The second property of the generalized coboundary matrices implies that the $h^{\text{th}}$ summand of the $l^{\text{th}}$ equation in (\ref{system}) reduces to be $r_{l+1,h} (m_{l+1,h}^i)^{x_i} (m_{l+1,h}^j)^{x_j}$, for $i$ and $j$ defining the (unique) two generalized coboundaries $\overline{M}_{\partial_i}$ and $\overline{M}_{\partial_j}$  sharing a negative entry in the position $(l+1,h)$. Notice that, eventually, one or even the two of these coboundaries $\partial _i,\partial _j$ might not be in {\bf B}.
\end{itemize}

\vspace{0.2cm}
Actually, the monomial $s_{l,h}(\text X)$ related to the mentioned $h^{\text{th}}$ summand of the $l^{\text{th}}$ equation in (\ref{system}) depends on whether the two, just one or none of the coboundaries  $\partial _i,\partial _j$ (precisely those whose related generalized coboundary matrices contribute a negative entry at position $(l+1,h)$) are in {\bf B}. More concretely,

\vspace{0.2cm}

\begin{itemize}
\item If $\partial_i,\partial_j, \in {\bf B}$, then
$$s_{l,h}(\text X):=r_{l+1,h}\,(1-2x_i)(1-2x_j).$$
\item If just $\partial_i \in {\bf B}$, then
$$s_{l,h}(\text X):=r_{l+1,h}\,(1-2x_i).$$
\item  If $\partial_i,\partial_j, \notin {\bf B}$, then
$$s_{l,h}(\text X):=r_{l+1,h}.$$
\end{itemize}
Let $S_l(\text X):=\sum_{j=1}^{4t}s_{l,j}(\text X)$ and let $\mathcal{H}^\rho_{G}$ be the set of solutions of (\ref{system}) of the form $\psi=\rho \prod _{i=1}^{k-m} \psi_i ^{x_i}$. The set $\mathcal{H}^\rho_{G}$ coincides with the set of solutions of the system of polynomial equations
\[\begin{cases}\begin{array}{rc}x_i(x_i-1)=0, & \,\, \text{ if } 1\leq i\leq k-m,\\
\quad S_l(\text X)=\sum_{j=1}^{4t} s_{l,j}(\text X)=0, & \,\,\text{ if }1\leq l\leq 4t-1.\end{array}\end{cases}\]

\vspace{0.2cm}

Similarly to Theorem \ref{thm1}, the next result holds.

\begin{thm} \label{thmbasesgeneral} The set $\mathcal{H}^\rho_{G}$ can be identified with the set of zeros of the following zero-dimensional ideal of $\mathbb{Q}[X]$.
$$J_G:=\langle\,x_i^2-x_i\colon\, i\in \{1,\ldots,k-m\}\,\rangle\, +\langle\,\sum_{h=1}^{4t} s_{l,h}(\text X)\colon\, l\in\{1,\ldots 4t-1\}\,\rangle.$$
Moreover, $|\mathcal{H}^\rho_{G}|= \mathrm{dim}_{\mathbb{Q}}(\mathbb{Q}[ X]/J_G)$.  \hfill $\Box$
\end{thm}

\vspace{0.5cm}

Observe in particular that, according to Lakshman and Lazard, the complexity of the computation of the reduced Gr\"obner decreases from $2^{O(t^2)}$ in Theorem \ref{thm1} to $2^{O(k-m)}$ in Theorem \ref{thmbasesgeneral}. In order to check the efficiency of this alternative, we have implemented in our library {\em hadamard.lib} a second procedure called {\tt CocCB($t,G,opt$)} that determines, depending on whether $opt=1$ or 2, the number or the explicit set of cocyclic Hadamard matrices developed over a given group $G$. The procedure has been tested in the computation of the number of cocyclic Hadamard matrices developed over the group $\mathbb{Z}_t \times \mathbb{Z}^2_2$ and the dihedral group $D_{4t}$ of order $4t$. Specifically, we have run the procedure in a system with an {\em Intel Core i7-2600, 3.4 GHz} and {\em Ubuntu}. Running times are exposed in Table \ref{table1}.

\begin{table}[h]
\begin{center}
\begin{tabular}{c|r|rr||r|rr|}
\ & \ & \multicolumn{2}{c|}{Running time in seconds} & \ & \multicolumn{2}{c}{Running time in seconds}\\
t & $|\mathcal{H}_{\mathbb{Z}_t\times \mathbb{Z}_2^2}|$ & {\em CocGM} & {\em CocCB} & $|\mathcal{H}_{D_{4t}}|$ & {\em CocGM} & {\em CocCB}\\ \hline
1 & 6 & 0 (0) & 0 (0) & 6 & 0 (0) & 0 (0)\\
3 & 24 & 129 (5718) & 0 (0) & 72 & - & 0 (0)\\
5 & 120 & - & 10 (120) & 1400 & - & 15 (-)\\
7 & - & - & - & 7488 & - & 68195 (-)\\
\end{tabular}
\end{center}

\caption{Running times related to {\em CocGM} and {\em CocCB}.}
\label{table1}
\end{table}

\vspace{0.2cm}

Actually, this procedure {\tt CocCB($t,G,opt$)} might be  improved if a deeper knowledge about the inner structure of cocyclic matrices over $G$ is known. In particular, progressing on the works of \cite{AGG15, AAFGGO16}, we have been able to design two specific procedures for looking for  $\mathbb{Z}_t \times \mathbb{Z}_2^2$-cocyclic Hadamard matrices and $D_{4t}$-cocyclic Hadamard matrices, so that larger cocyclic Hadamard matrices (up to $t \leq 31$) are obtained. The details are exposed in the next two subsections.

\subsection{The group $\mathbb{Z}_t\times \mathbb{Z}^2_2$}

Consider the group $G=\mathbb{Z}_t \times \mathbb{Z}_2^2$, $t>1$ odd, with ordering $$G=\{(0,0,0),(0,0,1),(0,1,0),(0,1,1),(1,0,0),\ldots,(t,1,1)\},$$ indexed as $\{1,\ldots,4t\}$. A basis ${\bf B}=\{\partial _2,\ldots, \partial _{4t-2},\beta_1,\beta_2,\beta_3\}$ for cocycles over $G$ is described by \cite{AAFR08, AAFR09}, and consists of $4t-3$ coboundaries and three representative cocycles. As usual, $\partial _i$ refers to the coboundary associated to the $i^{th}$-element in $G$. An explicit description of these cocycles was exposed by \cite{AAFR08}. Notice that all cocyclic Hadamard cocyclic matrices over $\mathbb{Z} _t \times \mathbb{Z}_2^2$ known so far use all the three representative cocycles $\beta_1,\beta_2$ and $\beta_3$ simultaneously (see the paper of \cite{BH95} for details). Thus, we assume
\[M_\rho=1_t \otimes  \left(\begin{array}{rrrr}1&1&1&1\\1&-1&1&-1 \\1&-1&-1&1\\1&1&-1&-1 \end{array}\right)\]
and hence, we restrict (\ref{0dime}) to the equations related to the $4t-3$ coboundaries, that is,
$$x_i(x_i-1)=0, \,\, \text{ for all } i\in\{1,\ldots,4t-3\}.$$

\vspace{0.2cm}

Let us point out that the system (\ref{system}) is equivalent to the one built up with the equations from the $4^{\text{th}}$ to the $(2t+1)^{\text{th}}$ (see the works of \cite{AAFR08, AGG15}). So, we have only $2t-2$ polynomials of the form $S_l(\text X)=s_{l,1}(\text X)+\ldots+s_{l,4t}(\text X)$. Before computing the monomials $s_{l,h}(\text X)$, we state the following lemma, which follows straightforwardly by inspection. In the sequel, $[n]_m$ denotes $n\mod m$ for short.

\begin{lemma}
Given the position $(s,c)$ with $5\leq s\leq 2t+2$ and $1\leq c\leq 4t$, the two generalized coboundary matrices with entries $-1$ at the position $(s,c)$ are $\overline{M}_{\partial_c}$ and $\overline{M}_{\partial_{j(s,c)}}$  where
\[j(s,c)  :=  1\,+\,4\displaystyle  \left[\left\lfloor \frac{s-1}{4}\right\rfloor+\left\lfloor\frac{c-1}{4}\right\rfloor\right]_t
  \,   +\,2\left[\left\lfloor\frac{[s-1]_4}{2}\right\rfloor \,+\,\left\lfloor\frac{ [c-1]_4}{2}\right\rfloor\right]_2
\,+\,
\left[s+c\right]_2.\] \hfill $\Box$
\end{lemma}

\vspace{0.2cm}

Taking into account this lemma and the basis ${\bf B}$ of cocycles, we compute the monomials
\[s_{l,h}(\text X):=r_{l+4,h}\,(1-2x_{h-1})^{\chi_{\bf B}(h)}\,(1-2x_{j(l+4,h)-1})^{\chi_{\bf B}(j(l+4,h))},\]
for $1\leq l\leq 2t-2$ and $1\leq h\leq 4t$, where
\[\chi_{\bf B}(i):=\left\{\begin{array}{cl}
1, & \text{if}\,\, \partial_i\in {\bf B},\\
0, & \text{otherwise}.
\end{array}\right.\]

\vspace{0.2cm}

Similarly to Theorem \ref{thm1}, the next result holds.

\begin{thm} \label{thm2} The set $\mathcal{H}^\rho_{\mathbb{Z}_t\times \mathbb{Z}^2_2}$ can be identified with the set of zeros of the following zero-dimensional ideal of $\mathbb{Q}[X]$.
$$J_{\mathbb{Z}_t\times \mathbb{Z}^2_2}:=\langle\,x_i^2-x_i\colon\, i\in \{1,\ldots,4t-3\}\,\rangle\, +\langle\,\sum_{h=1}^{4t} s_{l,h}(\text X)\colon\, l\in\{4,\ldots 2t+1\}\,\rangle.$$
Besides, $|\mathcal{H}^\rho_{\mathbb{Z}_t\times \mathbb{Z}^2_2}|= \mathrm{dim}_{\mathbb{Q}}(\mathbb{Q}[X]/J_{\mathbb{Z}_t\times \mathbb{Z}^2_2})$. \hfill $\Box$
\end{thm}

\vspace{0.2cm}

Actually, some additional assumptions, as those exposed by \cite{AGG15}, may be considered. Coboundaries $\partial _i$ on $\mathbb{Z}_t\times \mathbb{Z}^2_2$-cocyclic Hadamard matrices $M_\psi =  M_\rho \prod _{i\in I}  M_{\partial _i}$ are somehow symmetrically distributed, in the sense that the relation $4k+j \in I \Leftrightarrow 4t-4k+j \in I$, $1 \leq k \leq \frac{t-1}{2}$, $1 \leq j \leq 4$, holds. Furthermore, the number $c_k$ of coboundaries of each subset $\{ 4k+j\in I: \; 1 \leq j \leq 4\}$, for a fixed $2 \leq k \leq t$, satisfies $c_1 - c_k \equiv 1 \; \mbox{mod }2$. And the number $r_j$ of coboundaries of each subset $\{ 4k+j\in I: \; 1 \leq k \leq t\}$, for $1 \leq j \leq 4$, give rise to a tuple $dist=(r_1,r_2,r_3,r_4)$ (termed {\em distribution} by \cite{AGG15}) which certainly satisfies some additional properties.

\vspace{0.2cm}

Any $\mathbb{Z}_t\times \mathbb{Z}^2_2$-cocyclic matrix $M_\psi =  M_\rho \prod _{i\in I}  M_{\partial _i}$ may be uniquely identified as a $(4 \times t)$ binary matrix $D_\psi=(d_{jk})$ (termed {\em diagram} by \cite{AGG15}), such that $d_{jk}=1$ if and only if $4(k-1)+j \in I$. The conditions described above have a straightforward translation in terms of $D_\psi$. More concretely,

\vspace{0.2cm}

\begin{itemize}
\item column $i$ of $D_\psi$ is equal to column $t+2-i$, for $2 \leq i \leq t$.

\item the sum of column $2 \leq j\leq t$ of $D_\psi$ is of different parity of that of column 1 of $D_\psi$.

\item the sum of each row of $D_\psi$ gives the distribution $dist=(r_1,r_2,r_3,r_4)$.

\end{itemize}

\vspace{0.2cm}

Thus the method may be improved, as soon as the distribution $dist$ and the number of coboundaries per column $col=(c_2,\ldots,c_{\frac{t-1}{2}})$ are provided. We have implemented this method as a {\sc Singular} procedure called {\tt CocAH($t,col,dist,H$)}. Since exhaustive calculations are not feasible for $t\geq 11$, we have included in this procedure a new parameter $H=(x_2,\ldots,x_{2t+2})$, which determines which coboundaries are fixed ($x_i=1$ means $\partial _i$ is used, whereas $x_i=0$ implies $\partial _i$ is not used), and which of them are unknowns to be settled in the search (those corresponding to values $x_i=2$). This way, $\mathbb{Z}_t\times \mathbb{Z}^2_2$-cocyclic Hadamard matrices have been found up to $t \leq 31$, as Table \ref{table2} shows.

\begin{table}[h]
\begin{center}
\resizebox{\textwidth}{!}{
\begin{tabular}{c|rrrrr}
\ & \ & \ & \ & Running time & \ \\
t & col & dist & Initial coboundaries & in seconds & Final coboundaries\\ \hline
3 & 4 & 2,2,2,2 & - & 0 & 2,5,6,7,8,9,10\\
5 & 2,2 & 2,2,2,2 & - & 0 & 2,7,8,9,10,13,14\\
7 & 0,2,2 & 2,2,2,2 & - & 0 & 3,11,12,13,14,17,18,23,24\\
9 & 0,2,2,2  & 2,2,4,4 & - & 0 & 2,3,4,10,12,15,16,17,19,21,23,27,28,30,32\\
11 & 2,2,2,2,0 & 4,6,2,4 & - & 1 & 3,5,6,10,11,14,16,17,20,29,32,34,36,38,39,41,42\\
13 & 1,3,3,1,1,1 & 8,4,4,4 & - & 3 & 2,4,6,9,11,12,13,14,15,17,24,25,29,36,37,41,42,43,45,47,48,50\\
15 & 1,3,3,3,1,3,1 & 10,4,8,8 & - & 26 & 2,3,6,9,11,12,13,15,16,18,19,20,21,25,27,28,29,33,37,39,40,41,46,47,48,49,51,52,53,55,56,58\\
17 & 1,3,3,3,1,1,1,1 & 8,6,4,10 & 2 & 78 & 2,3,5,9,10,12,14,15,16,17,18,20,23,28,32,33,37,44,48,51,53,54,56,58,59,60,61,62,64,65\\
19 & 0,4,2,0,2,2,2,2,2 & 6,8,8,10 & 21 & 56 & 4,9,10,11,12,13,14,21,23,27,28,30,32,35,36,38,40,42,44,47,48,50,52,55,56,57,59,65,66,69,70,71,72\\
21 & 0,0,4,0,2,2,2,0,2,4 & 8,8,8,8 & - & 5 &
2,13,14,15,16,23,24,27,28,29,30,37,38,41,42,43,44,45,46,47,48,49,50,57,58,63,64,67,68,73,74,75,76\\
23 & 1,1,3,3,3,1,3,1,3,3,1 & 8,14,12,12 & 6,10,27,35 & 72 & 3,4,6,10,13,14,15,18,19,20,21,22,24,27,29,31,32,35,38,39,40,41,42,44,\\
\ & \ & \ & \ & \ &
48,52,53,54,56,58,59,60,63,65,67,68,71,73,74,76,78,79,80,81,82,83,86,90\\
25 & 2,2,2,2,4,2,2,2,2,0,2,2 & 14,14,8,12 & 9,12,14,16,18,20,25,27,30 & 158 & 4,5,6,9,12,14,16,18,20,21,22,23,24,25,27,30,32,33,34,38,39,45,47,49,52,\\
\ & \ & \ & \ & \ &
53,56,57,59,66,67,69,70,74,76,77,79,81,82,83,84,86,88,90,92,93,96,97,98\\
27 & 1,3,3,1,3,3,3,1,3,1,1,1,1 & 16,8,14,12 & 18,25,40,46,48,55 & 148 & 2,4,5,9,10,11,13,15,16,19,21,22,23,26,27,28,29,31,32,36,37,38,40,41,47,49,56,60,\\
\ & \ & \ & \ & \ &
61,67,69,73,74,76,80,81,83,84,86,87,88,89,90,91,95,97,99,100,101,102,103,105\\
29 & 1,1,3,1,1,1,3,1,3,3,3,3,3,1 & 14,12,18,12 & 5,9,18,23,25, & 36 & 2,3,5,9,14,15,16,18,23,25,29,31,32,35,38,39,40,41,43,44,46,47,48,49,51,52,53,54,55,58,62,\\
\ & \ & \ & \ & \ & 65,66,67,69,71,72,74,75,76,77,79,80,82,83,84,87,89,91,92,93,99,102,106,107,108,109,113\\
31 & 0,4,2,4,2,2,2,2,0,2,2,0,4,2,2 & 12,18,18,12 & 14,15,21,24,29 & 315 & 4,9,10,11,12,14,15,17,18,19,20,21,24,26,27,29,32,33,36,42,43,46,47,53,54,55,56,58,59,62,63,66,67,70,\\
\ & \ & \ & \ & \ & 71,73,74,75,76,82,83,86,87,93,96,97,100,102,103,105,108,109,110,111,112,114,115,117,118,119,120\\
\end{tabular}}
\end{center}
\caption{Auxiliary matrix method related to the group $\mathbb{Z}_t\times \mathbb{Z}^2_2$.}
\label{table2}
\end{table}

\subsection{The dihedral group $D_{4t}$}

Let $G$ be the dihedral group $D_{4t}=\langle a, b \,\colon\, a^{2t}=b^2=1,\,bab=a^{-1}\rangle$ with ordering
$$\{1,a,\ldots,a^{2t-1},b,ab,\ldots,a^{2t-1}b\},$$
indexed as $\{1,\ldots,4t\}$.  A basis ${\bf B}$ for cocycles over $G$ is explictly described by \cite{AAFR08, AAFR09}. For $t>2$, the basis consists of $4t-3$ coboundaries $\partial_k$ and three representative cocycles $\beta_i$, so that ${\bf B}=\{\partial _2,\ldots, \partial _{4t-2},\beta_1,\beta_2,\beta_3\}$. In the sequel we assume $t>2$. \cite{Fla97} observed that cocyclic Hadamard matrices over $D_{4t}$ mostly use $\beta_2\cdot \beta_3$ and do not use $\beta_1$. So, we assume $M_\rho=M_{\beta_2} \cdot M_{\beta_3} =\left( \begin{array}{cr} A&A\\B&-B  \end{array}\right)$, where
\[A=\left(
\begin{array}{crcr}
1&1&\cdots & 1\\ 1&&_{\cdot}\cdot^{\cdot}&-1\\ \vdots&_{\cdot}\cdot^{\cdot}&_{\cdot}\cdot^{\cdot}&\vdots\\ 1&-1&\cdots&-1\end{array} \right) \quad \mbox{and} \quad B=\left( \begin{array}{rrrr} 1&-1&\cdots & -1\\
\vdots&\ddots&\ddots&\vdots\\ 1&&\ddots&-1\\ 1&1&\cdots&1\end{array} \right).\] Therefore, in this case, (\ref{0dime}) is rewritten as
$$x_i(x_i-1)=0, \,\, \text{ for all } i\in\{1,\ldots,4t-3\}.$$

\vspace{0.2cm}

According to \cite{AAFR08}, the last $3t$ equations in the system (\ref{system}) are superfluous for $D_{4t}$-cocyclic Hadamard matrices. So, we have only $t-1$ polynomials of the form $S_l=s_{l,1}+\ldots+s_{l,4t}$. Before computing the monomials $s_{l,h}$, we state the following lemma, which follows straightforwardly by inspection.

\begin{lemma} The two generalized coboundary matrices with entries $-1$ in a position $(s,c)\in \{2,\ldots,t\}\times\{1,\ldots,4t\}$ are $\overline{M}_{\partial_c}$ and $\overline{M}_{\partial_{j(s,c)}}$, where

\vspace{0.2cm}

\begin{itemize}
\item If $1\leq c\leq 2t$, then
$$
j(s,c):=\left\{\begin{array}{cl}
2t, & \mbox{if $c+s-1=2t$},\\
c+s-1 \mod 2t, & \mbox{otherwise}.
\end{array}\right.$$
\item If $2t+1\leq c\leq 4t$, then
$$
j(s,c):=\left\{\begin{array}{cl}
c+s-1, & \mbox{if $c+s-1\leq 4t$},\\
2t + (c+s-1 \mod 2t), & \mbox{otherwise}.
\end{array}\right.$$ \hfill $\Box$
\end{itemize}
\end{lemma}

\vspace{0.2cm}

Taking into account this lemma and the basis of cocycles ${\bf B}$, we compute the monomials
$$s_{l,h}:=r_{l+1,h}\,(1-2x_{h-1})^{\chi_{\bf B}(h)}\,(1-2x_{j(l+1,h)-1})^{\chi_{\bf B}(j(l+1,h))},$$
for $1\leq l\leq t-1$ and $1\leq h\leq 4t$. Similarly to Theorem \ref{thm1}, the next result holds.

\begin{thm} \label{thm3} The set $\mathcal{H}^\rho_{D_{4t}}$ can be identified with the set of zeros of the following zero-dimensional ideal of $\mathbb{Q}[ X]$.
$$J_{D_{4t}}:=\langle\,x_i^2-x_i\colon\, i\in \{1,\ldots,4t-3\}\,\rangle\,+\,\langle\,\sum_{h=1}^{4t} s_{l,h}(\text X)\colon\, l\in\{1,\ldots t-1\}\,\rangle.$$
Besides, $|\mathcal{H}^\rho_{D_{4t}}|= \mathrm{dim}_{\mathbb{Q}}(\mathbb{Q}[ X]/J_{D_{4t}})$. \hfill $\Box$
\end{thm}

\vspace{0.2cm}

We have implemented this method as a {\sc Singular} procedure called {\tt CocDH($t,opt,H$)}. Once again, the parameter $opt$ indicates whether to compute the cardinality or the full set $\mathcal{H}^\rho_{D_{4t}}$ of $D_{4t}$-cocyclic matrices. And the auxiliary parameter $H=(x_2,\ldots,x_{4t-2})$ once again determines which coboundaries are fixed ($x_i=1$ means $\partial _i$ is used, whereas $x_i=0$ implies $\partial _i$ is not used), and which of them are unknowns to be settled in the search (those corresponding to values $x_i=2$). Similarly to the $\mathbb{Z}_t\times \mathbb{Z}^2_2$ case, an initial distribution of the coboundaries is fixed by considering the tuple $dist=(d_1,\ldots,d_t)$ such that

\vspace{0.2cm}

$$\begin{cases} d_1:= x_1+x_{2t-2} + x_{2t+1},\\
d_i:=x_i+x_{2t-i+1} + x_{2t+i} + x_{4t-i-1}, \text{ for all } i\in\{2,\ldots,t-1\},\\
d_t:=x_{2t-1}+x_{2t}.
\end{cases}$$

\vspace{0.2cm}

This way, $D_{4t}$-cocyclic Hadamard matrices have been found up to $t \leq 33$, as Table \ref{table3} shows. In the table, those initial coboundaries that are not permitted to be used are underlined.

\begin{table}[h]
\begin{center}
\resizebox{\textwidth}{!}{
\begin{tabular}{c|rrrr}
\ & \ & \ & Running time & \ \\
t & $dist$ & Initial coboundaries & in seconds & Final coboundaries\\ \hline
1 & 1 & - & 0 & 2\\
3 & 2,1,2 & - & 0 & 2,5,6,7,9\\
5 & 2,1,2,1,2 & - & 0 & 2,3,5,7,10,11,12,17\\
7 & 2,1,2,1,2,1,2 & - & 0 & 2,7,11,13,14,15,17,20,23,24,25\\
9 & 2,1,2,1,2,1,2,1,2 & - & 1 & 2,8,11,12,16,18,19,20,22,24,27,31,32,33\\
11 & 2,1,2,1,2,1,2,1,2,1,2 & - & 42 & 2,6,8,10,11,15,19,20,21,22,23,28,31,35,38,40,41\\
13 & 2,1,2,1,2,3,2,2,0,0,2,1,2 & - & 30 & 4,6,7,8,12,15,20,22,25,26,27,28,30,34,35,40,44,46,47,50\\
15 & 3,2,3,1,3,2,2,1,1,2,0,3,2,3,1 & \underline{16} & 108 & 2,3,4,5,6,7,8,11,13,15,18,21,24,25,29,30,32,33,34,36,39,41,44,45,46,\\
\ & \ & \ & \ & 47,48,53,57\\
17 & 1,2,2,3,3,2,3,1,4,4,2,1,2,2,2,3,2 & \underline{18},\underline{27} & 259 & 3,5,6,7,8,10,11,14,15,17,20,23,24,25,30,31,34,35,36,38,40,42,44,45,\\
\ & \ & \ & \ & 48,50,51,52,53,56,57,58,59,60,61,62,63,64,66\\
19 & 2,1,3,2,4,1,1,2,1,1,1,3,2,2,3,1,2,2,0 & \underline{4},7,8,12,\underline{13} & 141 & 2,6,7,8,12,15,18,23,26,29,30,33,34,35,40,41,42,43,44,51,52,53,54,55,\\
\ & \ & \ & \ & 57,58,59,61,63,64,66,68,71,73\\
21 & 3,2,1,2,1,1,2,1,1,1,3,1,3,1,2,4,2,0,1,3,2 & 6,9,\underline{12},13,75,78 & 76 & 2,6,8,9,13,14,16,17,22,26,29,31,35,41,42,43,44,45,46,47,53,54,56,59,\\
\ & \ & \ & \ & 60,62,63,64,67,68,69,70,73,75,78,80,82\\
23 & 1,3,2,2,4,3,2,0,3,2,3,2,1,2,2,3,1,3,2,1,2,1,1 & 14,29,\underline{40},\underline{63},67,70,\underline{74},\underline{81},\underline{83} & 213 & 2,3,6,7,10,12,13,14,15,16,17,19,20,22,25,27,28,29,30,35,36,37,39,41,\\
\ & \ & \ & \ & 46,49,50,51,52,53,54,56,58,61,62,65,67,70,76,80,82,86,87,88,89,90\\
25 & 1,3,4,1,2,1,3,1,2,3,1,1,2,1,1,2,0,1,2,2,2,2,3,3,1 & 9,19,\underline{23},\underline{24},\underline{26},50,63,65,66,89
& 272 & 3,4,8,9,10,14,17,19,20,21,22,25,27,28,29,30,34,37,40,43,44,46,47,48,\\
\ & \ & \ & \ & 50,52,54,56,58,60,61,63,65,66,70,73,74,75,76,77,89,90,95,97,98\\
27  & 0,3,4,2,3,3,2,2,2,4,0,1,1,2,2,2,2,1,2,3,2,2,3,1,2,1,2 & 19,\underline{21},\underline{24},25,28,29,\underline{48},68,96,\underline{103} & 117 & 4,6,7,8,9,11,16,18,19,20,22,25,28,29,31,33,34,37,38,44,45,49,50,51,52,\\
\ & \ & \ & \ & 54,55,57,58,59,60,61,62,63,64,65,68,69,70,75,77,78,83,85,86,88,89,92,\\
\ & \ & \ & \ & 94,96,98,102,105,106\\
29 & 2,2,2,2,2,3,3,3,2,3,1,3,4,1,2,4,2,2,1,2,2,2,2,1,3,1,3,3,1 & \underline{11},\underline{13},15,\underline{28},34,59,\underline{66},\underline{67},70,78,\underline{84},\underline{87},90,\underline{110}
& 203 & 2,4,6,7,8,9,10,14,15,16,17,21,22,26,29,30,31,33,34,37,40,42,43,45,46,\\
\ & \ & \ & \ & 48,50,51,52,53,54,55,56,57,59,61,63,65,68,69,70,71,72,75,76,78,79,\\
\ & \ & \ & \ & 81,82,86,88,89,90,91,93,94,98,99,100,103,104,106,108,109\\
31 & 1,2,2,2,4,4,1,1,1,3,1,2,0,2,2,2,4,2,3,3,3,1,3,3,0,2,2,2,3,1,1 & 2,9,\underline{20},\underline{21},\underline{25},\underline{28},\underline{33},40,51,\underline{52},63,72,\underline{86},\underline{93},94,\underline{103},104,117 & 625 & 2,5,6,7,9,11,15,16,18,22,24,30,35,36,38,39,40,41,42,43,45,46,51,56,57,\\
\ & \ & \ & \ & 59,63,65,68,69,72,73,75,79,80,81,82,83,84,87,89,91,92,94,95,96,97,\\
\ & \ & \ & \ & 100,101,104,105,106,107,109,110,112,114,117,118,119,120,121,122\\
33 & 2,1,4,1,2,2,3,2,3,1,1,1,2,3,0,2,3,3,1,2,3,2,4,2,1,2,1,0,2,2,1,0,0 & \underline{2},11,\underline{18},20,26,36,\underline{45},\underline{48},49,71,\underline{76},\underline{81},101,121,\underline{125},130 & 114 & 4,6,8,10,11,15,19,20,22,23,24,26,27,36,43,44,46,49,52,53,57,58,59,60,\\
\ & \ & \ & \ & 63,65,68,70,71,72,73,74,75,79,83,84,85,88,90,91,93,94,96,101,102,\\
\ & \ & \ & \ & 103,108,109,111,112,114,115,116,118,119,121,123,129,130
\end{tabular}}
\end{center}
\caption{Auxiliary matrix method related to the group $D_{4t}$.}
\label{table3}
\end{table}

\section{Conclusions}

By means of distinct techniques in algebraic geometry, this paper has been concerned with the computation of cocyclic Hadamard matrices over a fixed group $G$ of order $4t$, as the affine varieties of certain non-zero dimensional radical ideals. All the procedures that are described in the paper have been implemented in the open computer algebra system for polynomial computations {\sc Singular} and are included in the library {\em hadamard.lib}, which is available online at \texttt{http://personales.us.es/raufalgan/LS/hadamard.lib}.

\vspace{0.2cm}

Based on the classic cocyclic test of \cite{HdL95}, our first approach (Theorem \ref{thm1}),  has excessive complexity even for very small $t$. In order to improve the efficiency of this polynomial method, we have used recent results on the inner structure of a cocyclic matrix and we have defined a different ideal that also characterizes the set of $G$-cocyclic Hadamard matrices  (Theorem \ref{thmbasesgeneral}). Improved versions of this procedure ({\tt CocAH($t,col,dist,H$)} and {\tt CocDH($t,opt,H$)}, based on Theorems \ref{thm2} and \ref{thm3}, respectively) have been used to perform local searches for $\mathbb{Z}_t\times \mathbb{Z}_2^2$-cocyclic Hadamard matrices and $D_{4t}$-cocyclic Hadamard matrices, so that matrices of order up to $4t \leq 124$ have been found. To this end, an auxiliary data $H$ is needed to perform these local searches, for $t\geq 11$ (an exhaustive search is only feasible for $t<11$). More concretely, the list $H$ indicates which coboundaries are fixed (either used or not), and which of them are considered unknowns to be settled. A very interesting future work is trying to characterize if there exist some typo structures for $H$ such that the existence of cocyclic Hadamard matrices over either $\mathbb{Z}_t\times \mathbb{Z}_2^2$ or $D_{4t}$ is predicted.

\nocite{*}
\bibliographystyle{elsart-harv}
\bibliography{HadGrob}

\end{document}